\documentclass[12pt]{amsart}
\usepackage{amsmath, amsthm, amscd, amsfonts}

\theoremstyle{definition}

\theoremstyle{remark}

\numberwithin{equation}{section}
\newfont{\kh}{msbm10}

\begin{document}

\title{A Remark on the Mandl's Inequality}

\author{Mehdi Hassani}

\address{Mehdi Hassani, \newline Department of Mathematics, Institute for Advanced Studies in Basic
Sciences, P.O. Box 45195-1159, Zanjan, Iran}

\email{mmhassany@yahoo.com}



\begin{abstract}
In this note, first we refine Mandl's inequality. Then, we
consider the product $p_1p_2\cdots p_n$ and we refine some known
lower bounds for it, and we find some upper bounds for it by using
Mandl's inequality and its refinement and the AGM-Inequality.
\end{abstract}

\maketitle

%\tableofcontents

\section{Introduction}
As usual, let $p_n$ be the $n^{th}$ prime. The Mandl's inequality
\cite{dusart,rs} asserts that for every $n\geq 9$ we have
\begin{equation}\label{mandl}
\sum_{i=1}^np_i<\frac{n}{2}p_n.
\end{equation}
Considering the AGM Inequality \cite{rooin} and (\ref{mandl}), for
every $n\geq 9$ we obtain
$$
p_1p_2\cdots
p_n<\left(\frac{1}{n}\sum_{i=1}^np_i\right)^n<\left(\frac{p_n}{2}\right)^n.
$$
So, we have
\begin{equation}\label{pub}
p_1p_2\cdots p_n<\left(\frac{p_n}{2}\right)^n\hspace{10mm}(n\geq
9),
\end{equation}
where also holds true by computation for $5\leq n\leq 8$. In other
hand, one can get a trivial lower bound for that product using
Euclid's proof of infinity of primes; Letting $E_n=p_1p_2\cdots
p_n-1$ for every $n\geq 2$, it is clear that $p_n<E_n$. So, if
$p_n<E_n<p_{n+1}$ then $E_n$ should has a prime factor among
$p_1,p_2,\cdots, p_n$ which isn't possible. Thus $E_n\geq p_{n+1}$
and for every $n\geq 2$ we have
$$
p_1p_2\cdots p_n>p_{n+1}.
$$
In 1957 in \cite{rad-toep}, Bonse used elementary methods to show
that
$$
p_1p_2\cdots p_n>p_{n+1}^2\hspace{10mm}(n\geq 4),
$$
and
$$
p_1p_2\cdots p_n>p_{n+1}^3\hspace{10mm}(n\geq 5).
$$
In 1960 P\'{o}sa \cite{posa} proved that for every $k>1$ there
exists an $n_k$ such that for all $n\geq n_k$ we have
$$
p_1p_2\cdots p_n>p_{n+1}^k.
$$
In 1988, J. Sand\'{o}r found some inequalities of similar type;
For example he showed that for every $n\geq 24$ we have
$$
p_1p_2\cdots p_n>p_{n+5}^2+p_{[\frac{n}{2}]}^2.
$$
In 2000 Panaitopol \cite{panai} showed that in P\'{o}sa's result
we can get $n_k=2k$. More precisely, he proved that for every
$n\geq 2$ we have
$$
p_1p_2\cdots p_n>p_{n+1}^{n-\pi(n)},
$$
in which $\pi(x)=$ the number of primes $\leq x$. In this paper,
first we refine Mandl's inequality by showing
$\frac{n}{2}p_n-\sum_{i=1}^np_i>\frac{n^2}{14}$ for every $n\geq
10$. This refinement helps us to sharpen (\ref{pub}). Also, we
refine Panaitopol's result by proving
$$
p_1p_2\cdots p_n>p_{n+1}^{(1-\frac{1}{\log
n})(n-\pi(n))}\hspace{10mm}(n\geq 101).
$$
During proofs we will need some known results which we review them
briefly; we have the following known bounds for the function
$\pi(x)$, \cite{dusart}:
\begin{equation}\label{pi-lb}
\pi(x)\geq\frac{x}{\log x}\left(1+\frac{1}{\log
x}\right)\hspace{10mm}(x\geq 599),
\end{equation}
and
\begin{equation}\label{pi-ub}
\frac{x}{\log x}\Big(1+\frac{1.2762}{\log
x}\Big)\geq\pi(x)\hspace{10mm}(x\geq 2).
\end{equation}
For every $n\geq 53$, we have \cite{panai}
\begin{equation}\label{logpn+1}
\log p_{n+1}<\log n+\log\log n+\frac{\log\log n-0.4}{\log n}.
\end{equation}
Also, for every $n\geq 3$, we have \cite{robin}
\begin{equation}\label{tetapn}
\theta(p_n)>n\left(\log n+\log\log n-1+\frac{\log\log
n-2.1454}{\log n}\right),
\end{equation}
in which $\theta(x)=\sum_{p\leq x}\log p$. Specially,
$\theta(p_n)=\log(p_1p_2\cdots p_n)$ and this will act as a key
for approximating $p_1p_2\cdots p_n$. Finally, just for insisting
we note that base of all logarithms are $e$.

\section{Refinement of Mandl's Inequality}
To prove Mandl's inequality, Dusart (\cite{dusart}, page 50) uses
the following inequality
\begin{equation}\label{intpi}
\int_2^{p_n}\pi(t)dt\geq c+\frac{p_n^2}{2\log
p_n}\left(1+\frac{3}{2\log p_n}\right)\hspace{10mm}(n\geq 109),
\end{equation}
in which
$$
c=35995-3Li(599^2)+\frac{599^2}{\log 599}\approx -47.06746,
$$
and
$$
Li(x)=\lim_{\epsilon\rightarrow 0^+
}\left(\int_0^{1-\epsilon}\frac{dt}{\log
t}+\int_{1+\epsilon}^x\frac{dt}{\log t}\right),
$$
is logarithmic integral \cite{a-s}. Note that he has got
(\ref{intpi}) using (\ref{pi-lb}). Also, for using (\ref{intpi})
to prove Mandl's inequality we note that
$$
\int_2^{p_n}\pi(t)dt=\sum_{i=2}^n\big(p_i-p_{i-1}\big)(i-1)=
\sum_{i=2}^n\big(ip_i-(i-1)p_{i-1}\big)-\sum_{i=2}^np_i=np_n-\sum_{i=1}^np_i.
$$
Therefore, we have
\begin{equation}\label{npn-sumnpi}
np_n-\sum_{i=1}^np_i\geq c+\frac{p_n^2}{2\log
p_n}\left(1+\frac{3}{2\log p_n}\right)\hspace{10mm}(n\geq 109).
\end{equation}
Considering (\ref{pi-ub}) and (\ref{npn-sumnpi}), for every $n\geq
109$ we obtain
\begin{eqnarray*}
np_n-\sum_{i=1}^np_i&\geq& c+\frac{p_n^2}{2\log
p_n}\left(\frac{0.2238}{\log p_n}\right)+\frac{p_n^2}{2\log
p_n}\left(1+\frac{1.2762}{\log p_n}\right)\\&\geq&
c+0.1119\frac{p_n^2}{\log^2
p_n}+\frac{p_n}{2}\pi(p_n)=c+0.1119\frac{p_n^2}{\log^2
p_n}+\frac{n}{2}p_n.
\end{eqnarray*}
So, for every $n\geq 109$ we have
\begin{equation}\label{npn/2-sumnpi}
\frac{n}{2}p_n-\sum_{i=1}^np_i\geq c+0.1119\frac{p_n^2}{\log^2
p_n}.
\end{equation}
In other hand, we have the following bounds for $p_n$
(\cite{rs-62}, page 69)
$$
n\log n\leq p_n\leq n(\log n+\log\log n)\hspace{10mm}(n\geq 6).
$$
Combining these bounds with (\ref{npn/2-sumnpi}), for every $n\geq
109$ we yield that
$$
\frac{n}{2}p_n-\sum_{i=1}^np_i\geq c+\frac{0.1119(n\log
n)^2}{\log^2\big(n(\log n+\log\log n)\big)}.
$$
Now, for every $n\geq 21152$ we have $c+\frac{0.0119(n\log
n)^2}{\log^2(n(\log n+\log\log n))}>\frac{n^2}{14}$, and so we
obtain the following inequality for every $n\geq 21152$
$$
\frac{n}{2}p_n-\sum_{i=1}^np_i>\frac{n^2}{14}.
$$
By computation we observe that it holds also for $10\leq n\leq
21151$. Thus, we get the following refinement of Mandl's
inequality
\begin{equation}\label{refofmandl}
\sum_{i=1}^np_i<\frac{n}{2}p_n-\frac{n^2}{14}\hspace{10mm}(n\geq
10).
\end{equation}

\section{Approximation of the Product $p_1p_2\cdots p_n$}
Using (\ref{refofmandl}) and the AGM inequality we have
\begin{equation}\label{prub}
p_1p_2\cdots
p_n<\left(\frac{p_n}{2}-\frac{n}{14}\right)^n\hspace{10mm}(n\geq
10).
\end{equation}
Note that (\ref{prub}) holds also for $5\leq n\leq 9$. This yields
an upper bound for the product $p_1p_2\cdots p_n$. About lower
bound, as mentioned in introduction we show that
\begin{equation}\label{prlb}
p_1p_2\cdots p_n>p_{n+1}^{(1-\frac{1}{\log
n})(n-\pi(n))}\hspace{10mm}(n\geq 101).
\end{equation}
To prove this considering (\ref{pi-lb}), (\ref{logpn+1}) and
(\ref{tetapn}) it is enough to prove that
\begin{eqnarray*}
&~&\left(1-\frac{1-\frac{1}{\log n}}{\log n}-\frac{1-\frac{1}{\log
n}}{\log^2 n}\right)\left(\log n+\log\log n+\frac{\log\log
n-0.4}{\log n}\right)\\
&~&<\log n+\log\log n-1+\frac{\log\log n-2.1454}{\log
n}\hspace{10mm}(n\geq 599),
\end{eqnarray*}
which by putting $x=\log n$, is equivalent with:
$$
\frac{1.7454x^3+1.4x^2-0.4}{x^3+x^2-x-1}<\log
x\hspace{10mm}{x\geq\log 599},
$$
and trivially this holds true; because for $x\geq\log 599$ we have
$\frac{1.7454x^3+1.4x^2-0.4}{x^3+x^2-x-1}<1.7454$ and $1.85<\log
x$. Therefore, we yield (\ref{prlb}) for all $n\geq
599$. For $101\leq n\leq 598$ computation verifies it.\\
Finally, we use a refinement of the AGM inequality to get some
better bounds. In \cite{rooin-elem} Rooin shows that for any
non-negative real numbers $x_1\leq x_2\leq\cdots\leq x_n$ we have
\begin{equation}\label{refofagm}
A_n-G_n\geq\frac{1}{n}\sum_{k=2}^nA_{n-1}^{\frac{n-k}{n}}(x_n^{\frac{1}{n}}-A_{n-1}^{\frac{1}{n}})^k\geq
0,
\end{equation}
in which $nA_n=\sum_{i=1}^n x_i$ and $G_n^n=\prod_{i=1}^n x_i$.
For using this refinements we need Robin's inequality (see
\cite{dusart}, page 51) which gives a lower bound for the average
$\frac{1}{n}\sum_{i=1}^np_i$; for every $n\geq 2$ it asserts
\begin{equation}\label{robin}
p_{\lfloor\frac{n}{2}\rfloor}\leq\frac{1}{n}\sum_{i=1}^np_i.
\end{equation}
Applying (\ref{refofagm}) on $p_1<p_2<\cdots<p_n$ and using
relations (\ref{robin}) and (\ref{refofmandl}), for every $n\geq
10$ we obtain
\begin{equation}\label{prubref}
p_1p_2\cdots
p_n<\left\{\left(\frac{p_n}{2}-\frac{n}{14}\right)-\Omega(n)\right\}^n,
\end{equation}
in which
$$
\Omega(n)=\frac{1}{n}\sum_{k=2}^n
p_{\lfloor\frac{n-1}{2}\rfloor}^{\frac{n-k}{n}}\left\{p_n^{\frac{1}{n}}-\left(\frac{p_n}{2}-
\frac{n}{14}\right)^{\frac{1}{n}}\right\}^k>0.
$$
In fact, all members under summation are positive. So
$$
\Omega(n)>\frac{1}{n}\left\{p_n^{\frac{1}{n}}-\left(\frac{p_n}{2}-
\frac{n}{14}\right)^{\frac{1}{n}}\right\}^n>\frac{p_n}{2n}\left(2^{\frac{1}{n}}-1\right)^n.
$$
Using this bound for $\Omega(n)$ and considering (\ref{prubref}),
for every $n\geq 10$ we obtain
$$
p_1p_2\cdots
p_n<\left\{\frac{p_n}{2}\left(1-\frac{\left(2^{\frac{1}{n}}-1\right)^n}{n}
\right)-\frac{n}{14}\right\}^n.
$$

\section{On a Limit Concerning the Product $p_1p_2\cdots p_n$}

Some people believe that `` $e$ is The Master of All ''
\cite{mccartin}. There are some reasons, which one of them is the
the result $\lim_{n\rightarrow\infty}\prod_n^{\sharp}=e$ with
$(\prod_n^{\sharp})^{p_n}=p_1p_2\cdots p_n=e^{\theta(p_n)}$ (see
\cite{ruiz}). In fact, considering the Prime Number Theorem that
is $\prod_n^{\sharp}=e^{\frac{\theta(p_n)}{p_n}}=e+o(1)$, when
$n\rightarrow\infty$. In this section, we prove that
$\prod_n^{\sharp}=e+O(\frac{1}{\log^4(n\log n)})$, when
$n\rightarrow\infty$. It is known \cite{dusart} that for $x>1$, we
have
\begin{equation}\label{tb4}
|\theta(x)-x|<d\frac{x}{\log^4x},
\end{equation}
where $d=1717433$. Using this and
$\prod_n^{\sharp}=e^{\frac{\theta(p_n)}{p_n}}$, we obtain
$$
e^{-\frac{d}{\log^4p_n}}<\frac{\prod_n^{\sharp}}{e}<e^{\frac{d}{\log^4p_n}}\hspace{10mm}(n\geq
1).
$$
We have $e^{-\frac{d}{\log^4p_n}}>1-\frac{d}{\log^4p_n}$. Also,
for $p_n>5270747586811033$ a geometric approximation yields
$e^{\frac{d}{\log^4p_n}}<1+\frac{d}{\log^4p_n}+\frac{d^2}{2\log^4p_n(\log^4p_n-d)}$,
and so
$$
1-\frac{d}{\log^4p_n}<\frac{\prod_n^{\sharp}}{e}<1+\frac{d}{\log^4p_n}+\frac{d^2}{2\log^4p_n(\log^4p_n-d)}
\hspace{10mm}(p_n>5270747586811033).
$$
It is known \cite{rs-62} that $p_n>n\log n$ for every $n\geq 1$.
Using this, for every $p_n>5270747586811033$ we obtain
$$
1-\frac{d}{\log^4(n\log
n)}<\frac{\prod_n^{\sharp}}{e}<1+\frac{d}{\log^4(n\log
n)}+\frac{d^2}{2\log^4(n\log n)(\log^4(n\log n)-d)}.
$$
This describes $\lim_{n\rightarrow\infty}\prod_n^{\sharp}=e$
explicitly and also yields that
$\prod_n^{\sharp}=e+O(\frac{1}{\log^4(n\log n)})$, as we claimed.

% ------------------------------------------------------------------------

\end{document}